\documentclass[10pt,a4paper]{article}
\usepackage[cp1251]{inputenc}
\usepackage[english,russian]{babel}
\usepackage[T1]{fontenc}
\usepackage{amsmath}
\usepackage{amsthm}
\usepackage{amsfonts}
\usepackage{amssymb}
\usepackage{mathrsfs}
\usepackage[left=2.3cm, right=2.3cm, top=2.5cm, bottom=3cm]{geometry}

\title{Performances piecewise defined functions in analytic form, prime-counting function, $\xi$ sets}
\author{Oleh Kyrhan  \\ profugo.canis@gmail.com}

\begin{document}

\maketitle

\begin{abstract}
В статье рассматривается представление дискретных функций, определенных в аналитической форме без использования приближений, а именно функции Хевисайда, тождественной функции, дельта-функции Дирака и функции распределения простых чисел.

А также в статье введен и рассмотрен новый тип множеств ($\xi$-множества) посредством аналогии взятой из нахождения суммы ряда Гранди и других противоречий в математики и физике. С помощью $\xi$-множеств интерпретируется парадокс Рассела в системе аксиом наивной теории множеств.
\end{abstract}

\section{Введения}

В статье рассматривается вопрос представления дискретно определенных функций в аналитической форме без использования аппроксимации, а именно: функции Хевисайда, функция тождества, дельта-функция Дирака, функция распределения простых чисел и доказана теорема о представлении любой кусочно-заданной функции
%%%%%%%%%%%%%%%%%%%%%%%%%%%%%%%%%%%%%%%%%%%%%%%%%%%
\begin{equation}
t(x) = \left\{ {\begin{array}{*{20}{c}}
	{{t_0}(x),~x < {x_1}}~~~~~~~ \\ 
	{{t_1}(x),~{x_1} \le x < {x_2}} \\ 
	{...} ~~~~~~~~~~~~~~~~~~~~~~\\ 
	{{t_n}(x),~{x_n} \le x} ~~~~~~
	\end{array}} \right.
\end{equation}
%%%%%%%%%%%%%%%%%%%%%%%%%%%%%%%%%%%%%%%%%%%%%%%%%%%
где ${x_1} < {x_2} < ... < {x_n}$ точки изменения значения функции $t(x)$. Будут показаны функции из использованием несобственного и определенного интегралов. Например функция Хевисайда дискретные форы которой \cite{1,2,3}:

\begin{equation}
H_1(x) = \left\{ {\begin{array}{*{20}{c}}
			{0,~x < 0}\\
			{1,~x \ge 0}
		\end{array}} \right.
\end{equation}
%%%%%%%%%%%%%%%%%%%%%%%%%%%%%%%%%%
и
%%%%%%%%%%%%%%%%%%%%%%%%%%%%%%%%%%
\begin{equation}
H_2(x) = \left\{ {\begin{array}{*{20}{c}}
	{0,~ x < 0}\\
	{\frac{1}{2},~x = 0}\\
	{1,~x > 0}
\end{array}} \right.
\end{equation}
%%%%%%%%%%%%%%%%%%%%%%%%%%%%%%%%%%
имеют аналитические формы с использования аппроксимации \cite{1,2,3}:
%%%%%%%%%%%%%%%%%%%%%%%%%%%%%%%%%%
$$H(x) = \mathop {\lim }\limits_{k \to \infty } \left( {\frac{1}{2} + \frac{1}{\pi }arctg{\rm{ }}k \cdot x} \right)$$
%%%%%%%%%%%%%%%%%%%%%%%%%%%%%%%%%%
$$H(x) = \mathop {\lim }\limits_{k \to \infty } \frac{1}{{1 + {e^{ - 2k \cdot x}}}}$$
%%%%%%%%%%%%%%%%%%%%%%%%%%%%%%%%%%
$$H(x) = \mathop {\lim }\limits_{k \to \infty } \left( {\frac{1}{2} + \frac{1}{2}erf{\rm{ }}k \cdot x} \right)$$
%%%%%%%%%%%%%%%%%%%%%%%%%%%%%%%%%%
и интегральное представление с использованием аппроксимации и несобственного интеграла \cite{2,3}:
%%%%%%%%%%%%%%%%%%%%%%%%%%%%%%%%%%
$$H(x) = \mathop {\lim }\limits_{\varepsilon  \to 0 + } \frac{1}{{2\pi i}}\int\limits_{ - \infty }^\infty  {\frac{1}{{\tau  - i\varepsilon }}} {e^{ix\tau }}d\tau $$
%%%%%%%%%%%%%%%%%%%%%%%%%%%%%%%%%%

\section{Формулировка основных результатов (функции Хевисайда, функция тождества)}

Переведем в аналитическую форму без аппроксимаций функции (2), (3) и функцию тождества, с помощью которых переведем функцию распределения простых чисел в аналитическую форму без аппроксимаций.

Рассмотрим интеграл:
%%%%%%%%%%%%%%%%%%%%%%%%%%%%%%%%%%
$$\int\limits_0^\infty  {\frac{{{e^t}}}{{{{\left( {1 + {e^{ t}}} \right)}^2}}}} dt$$
%%%%%%%%%%%%%%%%%%%%%%%%%%%%%%%%%%
который равняется $\frac {1} {2} $. Теперь на основе этого интеграла построим функцию:
%%%%%%%%%%%%%%%%%%%%%%%%%%%%%%%%%%
\[f(x) = \int\limits_0^\infty  {\frac{{x{e^{xt}}}}{{{{\left( {1 + {e^{xt}}} \right)}^2}}}} dt\]

Найдем значение этой функции на числовой оси в соответствии со значениями переменой $x$.

%%%%%%%%%%%%%%%%%%%%%%%%%%%%%%%%%%
$$f(x) = \int\limits_0^\infty  {\frac{{x{e^{xt}}}}{{{{\left( {1 + {e^{xt}}} \right)}^2}}}} dt = \mathop {\lim }\limits_{t \to \infty } \left( { - \frac{1}{{1 + {e^{xt}}}}} \right) - \left( { - \frac{1}{{1 + {e^{x0}}}}} \right) = \frac{1}{2} - \mathop {\lim }\limits_{t \to \infty } \left( {\frac{1}{{1 + {e^{xt}}}}} \right)$$
%%%%%%%%%%%%%%%%%%%%%%%%%%%%%%%%%%
для $x>0$
%%%%%%%%%%%%%%%%%%%%%%%%%%%%%%%%%%
$$f(x > 0) = \frac{1}{2} - \mathop {\lim }\limits_{t \to \infty } \left( {\frac{1}{{1 + {e^{xt}}}}} \right) = \frac{1}{2} - 0 = \frac{1}{2}$$
%%%%%%%%%%%%%%%%%%%%%%%%%%%%%%%%%%
для $x<0$
%%%%%%%%%%%%%%%%%%%%%%%%%%%%%%%%%%
$$f(x < 0) = \frac{1}{2} - \mathop {\lim }\limits_{t \to \infty } \left( {\frac{1}{{1 + {e^{xt}}}}} \right) = \frac{1}{2} - 1 =  - \frac{1}{2}$$
%%%%%%%%%%%%%%%%%%%%%%%%%%%%%%%%%%
для $x=0$
%%%%%%%%%%%%%%%%%%%%%%%%%%%%%%%%%%
$$f(x = 0) = \frac{1}{2} - \mathop {\lim }\limits_{t \to \infty } \left( {\frac{1}{{1 + {e^{0t}}}}} \right) = \frac{1}{2} - \frac{1}{2} = 0$$

В дискретном определении функция $f(x)$ имеет следующую форму:
%%%%%%%%%%%%%%%%%%%%%%%%%%%%%%%%%%
$$f(x) = \left\{ {\begin{array}{*{20}{c}}
	{-\frac{1}{2},~x < 0}~~~ \\ 
	{0,~x = 0} \\ 
	{  \frac{1}{2},~x > 0} 
	\end{array}} \right.$$

Теперь прибавим $\frac{1}{2}$ к функции $f(x)$ и получим функцию которая совпадает из функциею (3)
%%%%%%%%%%%%%%%%%%%%%%%%%%%%%%%%%%
\begin{equation}
{H_2}(x) = f(x) + \frac{1}{2} = \int\limits_0^\infty  {\frac{{x{e^{xt}}}}{{{{\left( {1 + {e^{xt}}} \right)}^2}}}} dt + \frac{1}{2}
\end{equation}

Чтобы получить функцию (2) рассмотрим функцию
$$u(x) = \int\limits_0^\infty  {{x^2}\cdot{e^{ - t\cdot{x^2}}}} dt$$

Найдем значения функции $u(x)$ на числовой оси
%%%%%%%%%%%%%%%%%%%%%%%%%%%%%%%%%%
\begin{equation}
u(x) = \int\limits_0^\infty  {{x^2}\cdot{e^{ - t\cdot{x^2}}}} dt = \left( { - {e^{ - t\cdot{x^2}}}} \right)\left| {\begin{array}{*{20}{c}}
	\infty  \\ 
	0 
	\end{array}} \right. =  - \mathop {\lim }\limits_{t \to \infty } {e^{ - t\cdot{x^2}}} - \left( { - {e^{ - 0\cdot{x^2}}}} \right) =  - \mathop {\lim }\limits_{t \to \infty } {e^{ - t\cdot{x^2}}} + 1
\end{equation}
%%%%%%%%%%%%%%%%%%%%%%%%%%%%%%%%%%
у нас есть два варианты $ x=0 $ и $x \ne 0 $ так как в уравнения (5) входит квадрат $x$
%%%%%%%%%%%%%%%%%%%%%%%%%%%%%%%%%%
$$u(x = 0) =  - \mathop {\lim }\limits_{t \to \infty } {e^{ - t\cdot 0}} + 1 = -1 + 1 = 0$$
%%%%%%%%%%%%%%%%%%%%%%%%%%%%%%%%%%
$$u(x \ne 0,~{x^2}) =  - \mathop {\lim }\limits_{t \to \infty } {e^{ - t\cdot k}} + 1 = 0 + 1 = 1$$
%%%%%%%%%%%%%%%%%%%%%%%%%%%%%%%%%%

Дискретная форма функции $u(x)$ имеет следующий вид
%%%%%%%%%%%%%%%%%%%%%%%%%%%%%%%%%%
$$u(x) = \left\{ {\begin{array}{*{20}{c}}
	{0,~x = 0} \\ 
	{1,~x \ne 0} 
	\end{array}} \right.$$
%%%%%%%%%%%%%%%%%%%%%%%%%%%%%%%%%%
функция $ rt(x) = 1 - u(x)$ называется функция тождества.

Теперь с помощью функций (4) и $rt(x)$ построим функцию (2)
%%%%%%%%%%%%%%%%%%%%%%%%%%%%%%%%%%
\begin{multline} \nonumber
{H_1}(x) = {H_2}(x) + \frac{1}{2}rt(x) = \int\limits_0^\infty  {\frac{{x{e^{  xt}}}}{{{{\left( {1 + {e^{  xt}}} \right)}^2}}}} dt + \frac{1}{2} + \frac{1}{2}\left( {1 - \int\limits_0^\infty  {{x^2}{e^{ - t{x^2}}}} dt} \right) =\\
%%%%%%%%%%%%%%%%%%%%%%%%%%%%%%%%%%
= 1 - \frac{1}{2} \int\limits_0^\infty  {{x^2}{e^{ - t{x^2}}}} dt + \int\limits_0^\infty  {\frac{{x{e^{  xt}}}}{{{{\left( {1 + {e^{  xt}}} \right)}^2}}}} dt
\end{multline}
%%%%%%%%%%%%%%%%%%%%%%%%%%%%%%%%%%
$${H_1}(x) = 1 + \int\limits_0^\infty  {\left( {\frac{{x{e^{xt}}}}{{{{\left( {1 + {e^{xt}}} \right)}^2}}} - \frac{1}{2}{x^2} {e^{ - t{x^2}}}} \right)} dt$$

Теперь покажем аналитические формы без аппроксимаций и несобственного интеграла функции (2) и $rt(x)$

Рассмотрим функцию $c(x)$
%%%%%%%%%%%%%%%%%%%%%%%%%%%%%%%%%%
$$c(x) = \int\limits_0^{\frac{\pi }{2}} {\frac{{x{{\sec }^2}t{e^{ - x\tan t}}}}{{{{\left( {1 + {e^{ - x\tan t}}} \right)}^2}}}dt} $$

Вычислим значение функции $c(x) $ при разных $x$
%%%%%%%%%%%%%%%%%%%%%%%%%%%%%%%%%%
$$c(x) = \int\limits_0^{\frac{\pi }{2}} {\frac{{x{{\sec }^2}t{e^{ - x\tan t}}}}{{{{\left( {1 + {e^{ - x\tan t}}} \right)}^2}}}dt}  =  - \frac{1}{{1 + {{\text{e}}^{x\tan t}}}}\left| {\begin{array}{*{20}{c}}
	{\frac{\pi }{2}} \\
	\\ 
	0 
	\end{array}} \right. =  - \frac{1}{{1 + {{\text{e}}^{x\tan \frac{\pi }{2}}}}} - \left( { - \frac{1}{{1 + {{\text{e}}^{x\tan 0}}}}} \right)$$
%%%%%%%%%%%%%%%%%%%%%%%%%%%%%%%%%%
так как $\tan \frac{\pi }{2} = \infty$ и $\tan 0 = 0$ то функция $c(x)$ принимает те самые значения что и функция $f(x)$.
%%%%%%%%%%%%%%%%%%%%%%%%%%%%%%%%%%
$$c(x) = \left\{ {\begin{array}{*{20}{c}}
	{-\frac{1}{2},~x < 0}~~~ \\ 
	{0,~x = 0} \\ 
	{  \frac{1}{2},~x > 0} 
	\end{array}} \right.$$
%%%%%%%%%%%%%%%%%%%%%%%%%%%%%%%%%%

Заменим $f(x)$ на $c(x)$ в (4), тем самим получим аналитическую форму без аппроксимаций и несобственного интеграла функцию (3):
%%%%%%%%%%%%%%%%%%%%%%%%%%%%%%%%%%
\begin{eqnarray}
{H_2}(x) = \int\limits_0^{\frac{\pi }{2}} {\frac{{x{{\sec }^2}t{e^{ - x\tan t}}}}{{{{\left( {1 + {e^{ - x\tan t}}} \right)}^2}}}dt}  + \frac{1}{2}
\end{eqnarray}
Теперь найдем аналитическую форму без аппроксимаций и несобственного интеграла функции $rt(x)$, рассмотрим функцию 
%%%%%%%%%%%%%%%%%%%%%%%%%%%%%%%%%%
$$q(x)=\int\limits_{0}^{\frac{\pi }{2}}{{{\text{e}}^{-{{x}^{2}}\tan t}}{{x}^{2}}{{\sec }^{2}}t~dt}$$

Вычислим значение функции $q(x) $ при разных $x$
%%%%%%%%%%%%%%%%%%%%%%%%%%%%%%%%%%
\begin{multline} \nonumber
q(x) = \int\limits_0^{\frac{\pi }{2}} {{{\text{e}}^{ - {x^2}\tan t}}{x^2}{{\sec }^2}tdt}  = \left( { - {{\text{e}}^{ - {x^2}\tan t}}} \right)\left| {\begin{array}{*{20}{c}}
	{\frac{\pi }{2}} \\ 
	0 
	\end{array}} \right. = \mathop {\lim }\limits_{t \to \infty } \left( { - {{\text{e}}^{ - {x^2}\tan t}}} \right) - \left( { - {{\text{e}}^{ - {x^2}\tan 0}}} \right) =\\
%%%%%%%%%%%%%%%%%%%%%%%%%%%%%%%%%%
= \mathop {\lim }\limits_{t \to \infty } \left( { - {{\text{e}}^{ - {x^2}\tan t}}} \right) + 1
\end{multline}
%%%%%%%%%%%%%%%%%%%%%%%%%%%%%%%%%%
у нас есть два варианты $ x=0 $ и $x \ne 0 $ так как в функцию $q(x)$ входит квадрат $x$
%%%%%%%%%%%%%%%%%%%%%%%%%%%%%%%%%%
$$q(x = 0) = \mathop {\lim }\limits_{t \to \frac{\pi}{2} } \left( { - {{\text{e}}^{ - 0\tan t}}} \right) + 1 =  - 1 + 1=0$$
%%%%%%%%%%%%%%%%%%%%%%%%%%%%%%%%%%
$$q(x \ne 0,~k = {x^2}) = \mathop {\lim }\limits_{t \to \frac{\pi }{2}} \left( { - {{\text{e}}^{ - k\tan t}}} \right) + 1 = 0 + 1 = 1$$
%%%%%%%%%%%%%%%%%%%%%%%%%%%%%%%%%%
как видно функция $q(x)$ совпадает с $u(x)$. Функция тождества в таком случаи имеет вид
%%%%%%%%%%%%%%%%%%%%%%%%%%%%%%%%%%
\begin{eqnarray}
rt(x) = 1 - q(x) = 1 - \int\limits_0^{\frac{\pi }{2}} {{{\text{e}}^{ - {x^2}\tan t}}{x^2}{{\sec }^2}t~dt}
\end{eqnarray}

Теперь определим функцию $H_1(x)$ c помощью функций (6) и (7)
%%%%%%%%%%%%%%%%%%%%%%%%%%%%%%%%%%
\begin{multline} \nonumber
{H_1}(x) = {H_2}(x) + \frac{1}{2}rt(x) = \int\limits_0^{\frac{\pi }{2}} {\frac{{x{{\sec }^2}t{e^{ - x\tan t}}}}{{{{\left( {1 + {e^{ - x\tan t}}} \right)}^2}}}dt}  + \frac{1}{2} + \frac{1}{2}\left( {1 - \int\limits_0^{\frac{\pi }{2}} {{{\text{e}}^{ - {x^2}\tan t}}{x^2}{{\sec }^2}tdt} } \right) =\\
%%%%%%%%%%%%%%%%%%%%%%%%%%%%%%%%%%
 = 1 + \int\limits_0^{\frac{\pi }{2}} {\left( {\frac{{x{{\sec }^2}t{e^{ - x\tan t}}}}{{{{\left( {1 + {e^{ - x\tan t}}} \right)}^2}}} - \frac{1}{2}{{\text{e}}^{ - {x^2}\tan t}}{x^2}{{\sec }^2}t} \right)dt} 
\end{multline}
%%%%%%%%%%%%%%%%%%%%%%%%%%%%%%%%%%
\begin{eqnarray}
{H_1}(x) = 1 + \int\limits_0^{\frac{\pi }{2}} {\left( {\frac{{x{{\sec }^2}t{e^{ - x\tan t}}}}{{{{\left( {1 + {e^{ - x\tan t}}} \right)}^2}}} - \frac{1}{2}{{\text{e}}^{ - {x^2}\tan t}}{x^2}{{\sec }^2}t} \right)dt}
\end{eqnarray}

Теперь у нас есть все чтобы доказать теорему о представлении любой кусочно-заданной функции в аналитической форме без использования аппроксимации.

%%%%%%%%%%%%%%%%%%%%%%%%%%%%%%%%%%
\section{Кусочно-задание функции}

\newtheorem{Th}{Теорема}
\begin{Th}\label{theorem1}
Всякую кусочно-заданную функцию (1) можно представить в аналитической форме без аппроксимации, если функции $t_0(x),~t_1(x),~t_2(x)...~t_n(x)$ имеют аналитическую форму без аппроксимаций.
\end{Th}

\begin{proof}
Используем функцию (8). Возьмем два числа $a<b$ и построим функцию единичного импульса
%%%%%%%%%%%%%%%%%%%%%%%%%%%%%%%%%%
$$I(x,~a,~b) = H_1(x - a) - H_1(x - b)$$
%%%%%%%%%%%%%%%%%%%%%%%%%%%%%%%%%%
которая равняется 1 когда $a \le x < b$ и 0 в остальных случаях. Составим функцию (1) используя $I(x)$ и функции $t_0(x),~t_1(x),~t_2(x)...~t_n(x)$  
%%%%%%%%%%%%%%%%%%%%%%%%%%%%%%%%%%
$$t(x) = \left( {1 - {H_1}(x - {x_1})} \right){t_0}(x) + \sum\limits_{i = 2}^{n - 1} {I(x,{x_{i - 1}},{x_i})} {t_{i - 1}}(x) + {H_1}(x - {x_n}){t_n}(x)$$
%%%%%%%%%%%%%%%%%%%%%%%%%%%%%%%%%%
исходя из построения функции $t(x)$, теорема \ref{theorem1} доказана.
\end{proof}

%%%%%%%%%%%%%%%%%%%%%%%%%%%%%%%%%%
\section{Дельта-функция Дирака}

Представим дельта-функцию Дирака \cite{4, 5} в аналитической форме без аппроксимации через производную функции (8) по переменой $x$
%%%%%%%%%%%%%%%%%%%%%%%%%%%%%%%%%%
\begin{multline} \nonumber
\frac{d{{H}^{*}}(x)}{dx}=\frac{d\left( \int\limits_{0}^{\infty }{\left( {{e}^{-t}}-{{x}^{2}}\cdot {{e}^{-t\cdot {{x}^{2}}}}+\frac{x\cdot {{e}^{t\cdot x}}}{{{(1+{{e}^{t\cdot x}})}^{2}}} \right)}dt \right)}{dx}=\int\limits_{0}^{\infty }{\frac{d\left( {{e}^{-t}}-{{x}^{2}}\cdot {{e}^{-t\cdot {{x}^{2}}}}+\frac{x\cdot {{e}^{t\cdot }}}{{{(1+{{e}^{t\cdot }})}^{2}}} \right)}{dx}dt}=\\
%%%%%%%%%%%%%%%%%%%%%%%%%%%%%%%%%%	
 = \int\limits_0^\infty  {\left( {\frac{{{{\text{e}}^{t\cdot x}}}}{{{{(1 + {{\text{e}}^{t\cdot x}})}^2}}} - 2{{\text{e}}^{ - t\cdot{x^2}}}x} \right)} dt  - \int\limits_0^\infty  {\left( {\frac{{2{{\text{e}}^{2t\cdot x}}t\cdot x}}{{{{(1 + {{\text{e}}^{t\cdot x}})}^3}}}} \right)} dt + \int\limits_0^\infty  {\left( {\frac{{{{\text{e}}^{t\cdot x}}t\cdot x}}{{{{(1 + {{\text{e}}^{t\cdot x}})}^2}}} + 2{{\text{e}}^{ - t\cdot{x^2}}}t\cdot{x^3}} \right)} dt
\end{multline}
%%%%%%%%%%%%%%%%%%%%%%%%%%%%%%%%%%

\begin{equation}
	\frac{{d{H^*}(x)}}{{dx}} = \int\limits_0^\infty  {\left( {\frac{{{{\text{e}}^{t\cdot x}} + {{\text{e}}^{t\cdot x}}t\cdot x}}{{{{(1 + {{\text{e}}^{t\cdot x}})}^2}}} - 2{{\text{e}}^{ - t\cdot{x^2}}}x} \right)} dt + \int\limits_0^\infty  {\left( {2{{\text{e}}^{ - t\cdot{x^2}}}t\cdot{x^3} - \frac{{2{{\text{e}}^{2t\cdot x}}t\cdot x}}{{{{(1 + {{\text{e}}^{t\cdot x}})}^3}}}} \right)} dt   
\end{equation}
%%%%%%%%%%%%%%%%%%%%%%%%%%%%%%%%%%

Найдем значения функции (9)
\begin{multline} \nonumber
\frac{{d{H^*}(x)}}{{dx}} = \int\limits_0^\infty  {\left( {\frac{{{{\text{e}}^{t\cdot x}}}}{{{{(1 + {{\text{e}}^{t\cdot x}})}^2}}} - 2{{\text{e}}^{ - t\cdot{x^2}}}x} \right)} dt + \int\limits_0^\infty  {\left( { - \frac{{2{{\text{e}}^{2t\cdot x}}t\cdot x}}{{{{(1 + {{\text{e}}^{t\cdot x}})}^3}}} + \frac{{{{\text{e}}^{t\cdot x}}t\cdot x}}{{{{(1 + {{\text{e}}^{t\cdot x}})}^2}}}} \right)} dt +  \int\limits_0^\infty  {\left( {2{{\text{e}}^{ - t\cdot{x^2}}}t\cdot{x^3}} \right)} dt = \\
%%%%%%%%%%%%%%%%%%%%%%%%%%%%%%%%%%	
\left. { = \mathop {\lim }\limits_{A \to \infty } \left( { - \frac{t}{{{{(1 + {{\text{e}}^{t\cdot x}})}^2}}}} \right)} \right|\begin{array}{*{20}{c}}
A \\ 
0 
\end{array} -\left. { \mathop {\lim }\limits_{A \to \infty } \left( {2{{\text{e}}^{ - t\cdot{x^2}}}t\cdot x + \frac{t}{{1 + {{\text{e}}^{t\cdot x}}}}} \right)} \right|\begin{array}{*{20}{c}}
A \\ 
0 
\end{array} =\\
%%%%%%%%%%%%%%%%%%%%%%%%%%%%%%%%%%		
= \mathop {\lim }\limits_{t \to \infty } \left( { - \frac{t}{{{{(1 + {{\rm{e}}^{t \cdot x}})}^2}}}} \right) + \mathop {\lim }\limits_{t \to \infty } \left( { - 2{{\rm{e}}^{ - t \cdot {x^2}}}t \cdot x} \right) + \mathop {\lim }\limits_{t \to \infty } \left( {\frac{t}{{1 + {{\rm{e}}^{t \cdot x}}}}} \right) - 0 - 0 + 0 =\\ 
%%%%%%%%%%%%%%%%%%%%%%%%%%%%%%%%%%	
 = \mathop {\lim }\limits_{t \to \infty } \left( { - \frac{t}{{{{(1 + {{\rm{e}}^{t \cdot x}})}^2}}}} \right) + \mathop {\lim }\limits_{t \to \infty } \left( { - 2{{\rm{e}}^{ - t \cdot {x^2}}}t \cdot x} \right) + \mathop {\lim }\limits_{t \to \infty } \left( {\frac{t}{{1 + {{\rm{e}}^{t \cdot x}}}}} \right)=\\
%%%%%%%%%%%%%%%%%%%%%%%%%%%%%%%%%%
 = \mathop {\lim }\limits_{t \to \infty } \left( {\frac{{t{{\rm{e}}^{t \cdot x}}}}{{{{(1 + {{\rm{e}}^{t \cdot x}})}^2}}}} \right) + \mathop {\lim }\limits_{t \to \infty } \left( { - 2{{\rm{e}}^{ - t \cdot {x^2}}}t \cdot x} \right)
\end{multline}
%%%%%%%%%%%%%%%%%%%%%%%%%%%%%%%%%%
второй член при любом $x<\infty$  равняется нулю:
%%%%%%%%%%%%%%%%%%%%%%%%%%%%%%%%%%	
$$\mathop {\lim }\limits_{t \to \infty } \left( { - 2{{\rm{e}}^{ - t \cdot {x^2}}}t \cdot x} \right) = \mathop {\lim }\limits_{t \to \infty } \left( { - 2\frac{{t \cdot x}}{{{{\rm{e}}^{t \cdot {x^2}}}}}} \right) = 0.$$

%%%%%%%%%%%%%%%%%%%%%%%%%%%%%%%%%%	

Рассмотрим первый член, при $x > 0$, пускай $ \infty  > k > 0 $ і $ k = \left| x \right| $ 
%%%%%%%%%%%%%%%%%%%%%%%%%%%%%%%%%%	
\begin{multline} \nonumber
\mathop {\lim }\limits_{t \to \infty } \left( {\frac{{t{{\rm{e}}^{t \cdot k}}}}{{{{(1 + {{\rm{e}}^{t \cdot k}})}^2}}}} \right) = \mathop {\lim }\limits_{t \to \infty } \left( {\frac{{t{{\rm{e}}^{t \cdot k}}}}{{{{\rm{e}}^{2t \cdot k}} + 2{{\rm{e}}^{t \cdot k}} + 1}}} \right) = \left[ {\frac{\infty }{\infty }} \right] = \mathop {\lim }\limits_{t \to \infty } \left( {\frac{{{t^2}{{\rm{e}}^{t \cdot k}}}}{{{\rm{2}}t \cdot {{\rm{e}}^{2t \cdot k}} + 2t \cdot {{\rm{e}}^{t \cdot k}}}}} \right) = \\ 
%%%%%%%%%%%%%%%%%%%%%%%%%%%%%%%%%%
 = \mathop {\lim }\limits_{t \to \infty } \left( {\frac{{{t^2}{{\rm{e}}^{t \cdot k}}}}{{2t \cdot {{\rm{e}}^{t \cdot k}}({{\rm{e}}^{t \cdot k}} + 1)}}} \right) = \mathop {\lim }\limits_{t \to \infty } \left( {\frac{t}{{2({{\rm{e}}^{t \cdot k}} + 1)}}} \right) = 0.
\end{multline}
%%%%%%%%%%%%%%%%%%%%%%%%%%%%%%%%%%
при $x < 0$ , пускай $ \infty  > k > 0 $ і $ k = \left| x \right| $ 

\begin{multline} \nonumber
%%%%%%%%%%%%%%%%%%%%%%%%%%%%%%%%%%	
\mathop {\lim }\limits_{t \to \infty } \left( {\frac{{t{{\rm{e}}^{t \cdot x}}}}{{{{(1 + {{\rm{e}}^{t \cdot x}})}^2}}}} \right) = \mathop {\lim }\limits_{t \to \infty } \left( {\frac{{t{{\rm{e}}^{ - t \cdot k}}}}{{{{(1 + {{\rm{e}}^{ - t \cdot k}})}^2}}}} \right) = \\
%%%%%%%%%%%%%%%%%%%%%%%%%%%%%%%%%%	
=\mathop {\lim }\limits_{t \to \infty } \left( {\frac{t}{{{{\rm{e}}^{t \cdot k}}{{(1 + {{\rm{e}}^{ - t \cdot k}})}^2}}}} \right) = \mathop {\lim }\limits_{t \to \infty } \left( {\frac{t}{{{{\rm{e}}^{t \cdot k}}{{(1 + \frac{1}{{{{\rm{e}}^{t \cdot k}}}})}^2}}}} \right) = \\
%%%%%%%%%%%%%%%%%%%%%%%%%%%%%%%%%%	
= \mathop {\lim }\limits_{t \to \infty } \left( {\frac{t}{{{{\rm{e}}^{t \cdot k}}{{(1 + \frac{1}{{{{\rm{e}}^{t \cdot k}}}})}^2}}}} \right) = \mathop {\lim }\limits_{t \to \infty } \left( {\frac{t}{{{{\rm{e}}^{t \cdot k}}{{(1 + 0)}^2}}}} \right) = 0.
%%%%%%%%%%%%%%%%%%%%%%%%%%%%%%%%%%
\end{multline}
%%%%%%%%%%%%%%%%%%%%%%%%%%%%%%%%%%
при $x = 0$
%%%%%%%%%%%%%%%%%%%%%%%%%%%%%%%%%%	
$$\mathop {\lim }\limits_{t \to \infty } \left( {\frac{{t{{\rm{e}}^{t \cdot x}}}}{{{{(1 + {{\rm{e}}^{t \cdot x}})}^2}}}} \right) = \mathop {\lim }\limits_{t \to \infty } \left( {\frac{{t{{\rm{e}}^0}}}{{{{(1 + {{\rm{e}}^0})}^2}}}} \right) = \mathop {\lim }\limits_{t \to \infty } \left( {\frac{t}{4}} \right) = \infty. $$

Как видно уравнения (9) имеет значения дельта-функции Дирака т.е. $\frac{{d{H^*}(x)}}{{dx}}~=~{\delta ^*}(x)$.

\section{Функция распределения простых чисел}

Чтобы получить функцию распределения простых чисел \cite{6}, нужно построить функцию количества делителей ${\sigma _0}(n)$ числа $n$
на основе которой строится функция идентификации простых чисел.

Переведем функцию ${\sigma _0}(n)$ количества делителей числа $n$ в аналитическую форму без аппроксимации. Используем свойство функции $\sin (x)$, если $x$ целое то $\sin (x) = 0$.

Если $i$  делит число  $n$ то $\sin \left( {\pi \frac{n}{i}} \right) = 0$  в противном случаи $\sin \left( {\pi \frac{n}{i}} \right) \ne 0$. Используем функцию (7), тогда следующая функция 
$rt\left( {\sin \left( {\pi \frac{n}{i}} \right)} \right)$ равняется 1 если $i$  делит $n$ и 0 в противном случаи. Теперь построим функцию  ${\sigma _0}(n)$ которая суммирует $rt\left( {\sin \left( {\pi \frac{n}{i}} \right)} \right)$ по всем $i$ 
%%%%%%%%%%%%%%%%%%%%%%%%%%%%%%%%%%
$${\sigma _0}(n) = \sum\limits_{i = 1}^\infty  {rt\left( {\sin \left( {\pi \frac{n}{i}} \right)} \right)} $$
 
Теперь построим функцию идентификации простых чисел на основе функции ${\sigma _0}(n)$ и $rt(x)$, так как у простого числа всего два делителя, 1 и оно само, то
%%%%%%%%%%%%%%%%%%%%%%%%%%%%%%%%%%
$$fes(n) = rt\left( {{\sigma _0}(n) - 2} \right) = rt\left( {\sum\limits_{i = 1}^\infty  {rt\left( {\sin \left( {\pi \frac{n}{i}} \right)} \right)}  - 2} \right)$$
%%%%%%%%%%%%%%%%%%%%%%%%%%%%%%%%%%
дискретная форма которого будет иметь следующий вид
%%%%%%%%%%%%%%%%%%%%%%%%%%%%%%%%%%
$$fes(n) = \left\{ {\begin{array}{*{20}{c}}
	{1,~n~\text{простое}}~~~ \\ 
	{0,~n~\text{составное}} 
	\end{array}} \right.$$
%%%%%%%%%%%%%%%%%%%%%%%%%%%%%%%%%%

Теперь имея функцию идентификации простых чисел и функцию Хевисайда $H_1(x)$ построим функцию распределения простых чисел в аналитической форме без использования аппроксимации.
%%%%%%%%%%%%%%%%%%%%%%%%%%%%%%%%%%
$$\pi (x) = \sum\limits_{i = 1}^\infty  {fes(i)} {H_1}(x - i) = \sum\limits_{i = 1}^\infty  {rt\left( {\sum\limits_{j = 1}^\infty  {rt\left( {\sin \left( {\pi \frac{i}{j}} \right)} \right)}  - 2} \right)} {H_1}(x - i)$$
%%%%%%%%%%%%%%%%%%%%%%%%%%%%%%%%%%

Общий вид функции распределения простых чисел:

{\tiny \begin{multline} \nonumber
\pi (x) =   \\
%%%%%%%%%%%%%%%%%%%%%%%%%%%%%%%%%%
= \sum\limits_{i = 1}^\infty  {\left( {1 - \int\limits_0^{\frac{\pi }{2}} {{{\text{e}}^{ - {{\left( {\sum\limits_{j = 1}^\infty  {\left( {1 - \int\limits_0^{\frac{\pi }{2}} {{{\text{e}}^{ - \sin {{\left( {\pi \frac{n}{j}} \right)}^2}\tan z}}\sin {{\left( {\pi \frac{n}{j}} \right)}^2}{{\sec }^2}zdz} } \right)}  - 2} \right)}^2}\tan v}}{{\left( {\sum\limits_{j = 1}^\infty  {\left( {1 - \int\limits_0^{\frac{\pi }{2}} {{{\text{e}}^{ - \sin {{\left( {\pi \frac{n}{j}} \right)}^2}\tan z}}\sin {{\left( {\pi \frac{n}{j}} \right)}^2}{{\sec }^2}zdz} } \right)}  - 2} \right)}^2}{{\sec }^2}vdv} } \right)}  \cdot  \\ 
%%%%%%%%%%%%%%%%%%%%%%%%%%%%%%%%%%
 \cdot  \left( {1 + \int\limits_0^{\frac{\pi }{2}} {\left( {\frac{{\left( {x - i} \right){{\sec }^2}t{e^{ - \left( {x - i} \right)\tan t}}}}{{{{\left( {1 + {e^{ - \left( {x - i} \right)\tan t}}} \right)}^2}}} - \frac{1}{2}{{\text{e}}^{ - {{\left( {x - i} \right)}^2}\tan t}}{{\left( {x - i} \right)}^2}{{\sec }^2}t} \right)dt} } \right)
\end{multline} 
%%%%%%%%%%%%%%%%%%%%%%%%%%%%%%%%%%
}

\noindent
Все вычисления были проверены в \textit {wolfram mathematica}.

% % % % % % % % % % % % % % % % % % % % % % % % % % % % % % %
% % % % % % % % % % % % % % % % % % % % % % % % % % % % % % %
% % % % % % % % % % % % % % % % % % % % % % % % % % % % % % %
% % % % % % % % % % % % % % % % % % % % % % % % % % % % % % %
% % % % % % % % % % % % % % % % % % % % % % % % % % % % % % %\\
% % % % % % % % % % % % % % % % % % % % % % % % % % % % % % %

\section{Введения 2} 

\begin{flushright}
\begin{minipage}[t]{100mm}
{\footnotesize Show me the infinity and I will prove the inconsistency of the Universe, let me infinity and I create my Universe.}
\begin{flushright}
{\footnotesize - \textit{Ron Swanson} -}
\end{flushright}
\end{minipage}
\end{flushright}

\noindent
\textbf{Противоречия в математике}

Некоторые противоречия\footnote{~В работе противоречия рассматриваются не как недостаток логики или неполнота теории, а как нечто что не вписывается в рамки логики и имеет право на существование.} используются на практике (или используются утверждения приводящие к противоречиям) в математике и физике, например: 

\noindent
$\bullet$ Использование мнимой единицы \cite{7} не нужно перечислять, так как без нее не было некоторых разделов математики и физики, но не только ее определения не вписывается в наше понимание, вона еще и приводит к противоречию:

$i = \sqrt { - 1}  = \sqrt {\frac{{ - 1}}{1}}  = \sqrt {\frac{1}{{ - 1}}}  = \frac{1}{{\sqrt { - 1} }} = \frac{1}{i}$ откуда следует $- 1 = {i^2} = ii = \frac{i}{i} = 1$.

\noindent
$\bullet$ Конечная сумма всех натуральных чисел \cite{3}
%%%%%%%%%%%%%%%%%%%%%%%%%%%%%%%%%%%%%%%%
$$\sum\limits_{n = 1}^\infty  n  = 1 + 2 + 3 + 4 + 5 + 6 + ...$$

\noindent
используется в объяснении эффекта Казимира и в теории струн. Существует множество способов найти сумму всех натуральных чисел, рассмотрим один из них: 
%%%%%%%%%%%%%%%%%%%%%%%%%%%%%%%%%%%%%%%%
$$~~~~c = 1 + 2 + 3 + 4 + 5 + 6 + ...$$ 
$$~~~4c =~~~~~4 + ~~~~~8 + ~~~~~12 + ...$$
$$-3c = 1 - 2 + 3 - 4 + 5 - 6 + ...$$
%%%%%%%%%%%%%%%%%%%%%%%%%%%%%%%%%%%%%%%%

\noindent
ряд $1 - 2 + 3 - 4 + 5 - 6 + ...$ является разложения в степенной ряд функции $1/{\left( {1 + x} \right)^2}$ при $x$, равном 1. Соответственно
%%%%%%%%%%%%%%%%%%%%%%%%%%%%%%%%%%%%%%%%
$$-3c = 1 - 2 + 3 - 4 + 5 - 6 + ... = 1/\left( {1 + 1} \right)^2 = \frac{1}{4}$$ 

\noindent
поделив обе части на $-3$ получаем $c=- \frac{1}{12}$.

Следствиям такого суммирования, это нахождения сумм следующих бесконечных рядов:
%%%%%%%%%%%%%%%%%%%%%%%%%%%%%%%%%%%%%%%%
$$\sum\limits_{n = 1}^\infty \frac{n}{n}  = 1 + 1 + 1 + 1 + 1 + ... =  - \frac{1}{2};$$
%%%%%%%%%%%%%%%%%%%%%%%%%%%%%%%%%%%%%%%%
$$\sum\limits_{n = 0}^\infty  {{n^2}}  = 1 + 4 + 9 + 16 + ... = 0.$$

Рассмотрим еще несколько противоречий которые имеют место в математике.

\noindent
$\bullet$ Теорема Римана об условно сходящихся рядах \cite{10} которая гласит что:

\textit{Пусть ряд $\mathbf{A}$ сходится условно, тогда для любого числа $\mathbf{S}\in\mathbb{R} \cup \{ \infty \}$ можно так поменять порядок суммирования, что сумма нового ряда будет равна $\mathbf{S}$.}

\noindent
$\bullet$ Ряд Гранди \cite{8,9}~--- это бесконечный ряд 
%%%%%%%%%%%%%%%%%%%%%%%%%%%%%%%%%%%%%%%%
$$1-1+1-1+1-1+1-1+... \text{~или~} 
\sum\limits_{n = 0}^\infty  {{{( - 1)}^n}}$$. 

Один из очевидных методов нахождения суммы ряда, это воспринимать его как телескопический ряд и попарно сгруппировать члены:
$(1-1)+(1-1)+(1-1)+(1-1)+...=0+0+0+0+...=0.$ С другой стороны, похожим способом можно получить другой ответ:
$1+(-1+1)+(-1+1)+(-1+1)+...=1+0+0+0+...=1.$

Таким образом, различной расстановкой скобок в ряде Гранди, можно получить в качестве суммы 0 или 1.
Если считать ряд Гранди расходящейся геометрической прогрессией, то, используя те же методы что и при работе со сходящимися геометрическими прогрессиями, можно получить третье значение, $1/2$:
%%%%%%%%%%%%%%%%%%%%%%%%%%%%%%%%%%%%%%%%
$$\sum\limits_{n = 0}^\infty  {{{( - 1)}^n}}  = 1 - 1 + 1 - 1 + ...$$
%%%%%%%%%%%%%%%%%%%%%%%%%%%%%%%%%%%%%%%%
$$1 - \sum\limits_{n = 0}^\infty  {{{( - 1)}^n}}  = 1 - (1 + 1 - 1 + ...) = 1 - 1 + 1 - 1 + ... = \sum\limits_{n = 0}^\infty  {{{( - 1)}^n}}$$
%%%%%%%%%%%%%%%%%%%%%%%%%%%%%%%%%%%%%%%%
$$1 - \sum\limits_{n = 0}^\infty  {{{( - 1)}^n}}  = \sum\limits_{n = 0}^\infty  {{{( - 1)}^n}}$$
%%%%%%%%%%%%%%%%%%%%%%%%%%%%%%%%%%%%%%%%
$$\sum\limits_{n = 0}^\infty  {{{( - 1)}^n}}  = \frac{1}{2}$$

\noindent
из этого можно прийти к двум выводам: Ряд $1 - 1 + 1 - 1 + ...$ не имеет суммы или его сумма должна быть равна $1/2$.

Перенесем (по аналогии) рассуждения про сумму ряда Гранди на множества.

%%%%%%%%%%%%%%%%%%%%%%%%%%%%%%%%%%%%%%%%
%%%%%%%%%%%%%%%%%%%%%%%%%%%%%%%%%%%%%%%%
\textbf{$\xi$-парадокс}

Для рассмотрения $\xi$-парадокса нам понадобится вспомогательная теорема об операциях над множествами. В тексте  пустое множество обозначается символом $\theta $.

\begin{Th}\label{theorem2}
	Для двух множеств $\mathbf A$ и $\mathbf B$, $\mathbf A \ne \mathbf B$ справедливо равенство 
	%%%%%%%%%%%%%%%%%%%%%%%%%%%%%%%%%%%%%%%%
	$$(\mathbf A \cap \mathbf B) \cup (\mathbf A \cap \mathbf B) \cup (\mathbf A \cap \mathbf B) \cup (\mathbf A \cap \mathbf B) \cup ... = \mathbf A \cap (\mathbf B \cup \mathbf A) \cap (\mathbf B \cup \mathbf A) \cap (\mathbf B \cup \mathbf A) \cap (\mathbf B \cup ...$$
\end{Th}

\begin{proof}
	Для двух множеств $\mathbf A$ и двух множеств $\mathbf B$, операции объединения и пересечения ассоциативная:
	%%%%%%%%%%%%%%%%%%%%%%%%%%%%%%%%%%%%%%%%
	$$\mathbf A \cap \mathbf B  = \mathbf A \cap \mathbf B  \cup \mathbf A \cap \mathbf B  = (\mathbf A \cap \mathbf B ) \cup (\mathbf A \cap \mathbf B) = \mathbf A \cap (\mathbf B  \cup \mathbf A) \cap \mathbf B $$
	
	Для трех множеств $\mathbf A$ и двух множеств  $\mathbf B$, операции объединения и пересечения ассоциативная: 
	%%%%%%%%%%%%%%%%%%%%%%%%%%%%%%%%%%%%%%%%
	$$\mathbf {A = A \cap B  \cup A \cap B  \cup A = (A \cap B ) \cup (A \cap B) \cup A = A \cap (B \cup A) \cap (B \cup A)}$$
	
	Для производного количества элементов ассоциативность поочередного использования операций объединения и пересечения доказывается по индукции. 
\end{proof}

\begin{Th}[$\xi$ парадокс] \label{theorem4}
	Для произвольного множества $\mathbf G \ne \theta$ справедливо утверждения 
	
	$\mathbf {G \ne \theta  \Leftrightarrow G = \theta }$.
\end{Th}

\begin{proof}
	Пустое множество представим, как бесконечное объединения пустых множеств:
	%%%%%%%%%%%%%%%%%%%%%%%%%%%%%%%%%%%%%%%%
	$$\mathbf { \theta  = \theta  \cup \theta  \cup \theta  \cup \theta  \cup \theta  \cup \theta  \cup ...}$$
	
	каждое из них, представим как пересечения произвольного не пустого множества $\mathbf G \ne \theta$ из пустым $\mathbf {\theta = G \cap \theta}$, из этого получим
	%%%%%%%%%%%%%%%%%%%%%%%%%%%%%%%%%%%%%%%%
	$$\mathbf { \theta = (G \cap \theta ) \cup (G \cap \theta ) \cup (G \cap \theta ) \cup (G \cap \theta ) \cup (G \cap \theta ) \cup (G \cap \theta ) \cup ...}$$
	%%%%%%%%%%%%%%%%%%%%%%%%%%%%%%%%%%%%%%%%
	
	\noindent
	согласно теореме \ref{theorem2} поменяем порядок поочередного применения операций пересечения и объединения
	%%%%%%%%%%%%%%%%%%%%%%%%%%%%%%%%%%%%%%%%
	$$\mathbf { \theta = G \cap (\theta  \cup G) \cap (\theta  \cup G) \cap (\theta  \cup G) \cap (\theta  \cup G) \cap (\theta  \cup G) \cap (\theta  \cup ... =}$$
	%%%%%%%%%%%%%%%%%%%%%%%%%%%%%%%%%%%%%%%%
	$$\mathbf {= G \cap G \cap G \cap G \cap G \cap G \cap ... = G. }$$
		
	Обратное утверждения доказывается аналогично, нужно рассмотреть эту процедуру в обратном порядке.
\end{proof}
%%%%%%%%%%%%%%%%%%%%%%%%%%%%%%%%%%%%%%%%

Теорема \ref{theorem4} утверждает что ``всякое непустое множество является пустым, и наоборот'' это противоречия, но по аналогии из рядом Гранди мы не станем отвергать теорему \ref{theorem4}, мы сделаем предположения о возможности существования множества которое одновременно может равняться двум или более множествам.

%%%%%%%%%%%%%%%%%%%%%%%%%%%%%%%%%%%%%%%%
%%%%%%%%%%%%%%%%%%%%%%%%%%%%%%%%%%%%%%%%
\section{$\xi$-множество}

\textbf{Определения $\xi$-множества.} Из теоремы \ref{theorem4} следует что ряд
%%%%%%%%%%%%%%%%%%%%%%%%%%%%%%%%%%%%%%%%	
\begin{eqnarray}\label{xi}
\mathbf {G \cap \theta  \cup G \cap \theta  \cup G \cap \theta  \cup G \cap \theta  \cup G \cap \theta  \cup G \cap \theta  \cup ...}
\end{eqnarray}
%%%%%%%%%%%%%%%%%%%%%%%%%%%%%%%%%%%%%%%%
равняется сразу двум множествам $\mathbf G$ и $\theta$ одновременно, как и ряд Гранди $0$ и $1$, предположим что ряд (\ref{xi}) равняется некоторому множеству назовем его $\xi$-множество.

Расширим ряд (\ref{xi}) на случай произвольных двух множеств, построим $\xi$-множество класса 2 которое одновременно равняется двум не пустым множествам. Возьмем два произвольных множества $\mathbf A \ne \theta $ и $\mathbf B \ne \theta $ из условием что 
%%%%%%%%%%%%%%%%%%%%%%%%%%%%%%%%%%%%%%%%
$$\mathbf {A \cap B = F \ne \theta }\text{~и~}\mathbf {A \cup B = D \ne \theta }$$

\noindent
на их основе построим ряд 
%%%%%%%%%%%%%%%%%%%%%%%%%%%%%%%%%%%%%%%%
$$\mathbf {A \cap B \cup A \cap B \cup A \cap B \cup A \cap B \cup A \cap B \cup ...}$$

\noindent
найдем какому множеству равняется этот ряд. Расставим дужки, и это даст нам следующий результат 
%%%%%%%%%%%%%%%%%%%%%%%%%%%%%%%%%%%%%%%%
$$\mathbf {A \cap (B \cup A) \cap (B \cup A) \cap (B \cup A) \cap (B \cup A) \cap (B \cup ... =}$$
%%%%%%%%%%%%%%%%%%%%%%%%%%%%%%%%%%%%%%%%
$$\mathbf {= A \cap D \cap D \cap D \cap D \cap ... = A}$$

\noindent
теперь согласно теореме \ref{theorem2} переставим дужки и получим 
%%%%%%%%%%%%%%%%%%%%%%%%%%%%%%%%%%%%%%%%
$$\mathbf {(A \cap B) \cup (A \cap B) \cup (A \cap B) \cup (A \cap B) \cup (A \cap B) \cup ...= }$$
%%%%%%%%%%%%%%%%%%%%%%%%%%%%%%%%%%%%%%%%
$$\mathbf {= F \cup F \cup F \cup F \cup F \cup ... = F}$$

\noindent
как видно этот ряд одновременно равняется двум множествам $\mathbf A$ и $\mathbf F$.

Из этого можно сделать определение $\xi$ множества

\newtheorem{Def}{Определение}
\begin{Def}
$\xi$ множество класса n это множество которое одновременно равняется нескольким множествам, класс $\xi$ множества это количество множеств которым равняется это $\xi$ множество.

Обозначим $\xi$-множество явно $\mathbf {G_1 || G_2 || G_3 ||...||G_n}$ (явно показывает каким именно множествам равняется $\xi$-множество) и не явно $\mathbf { \widetilde G^n}$.

Если класс $\xi$-множества бесконечный, то он обозначается алефом того множества $\aleph$ множеств котором оно равняется одновременно. 

Обычное множество это $\xi$-множество произвольного класса которое равняется одному и то муже множеству $\mathbf {A = \widetilde A = A||A||A||...||A||...}$
\end{Def}

\noindent
\textbf{Функции образования $\xi$-множеств, $\xi$-функции} 
\begin{Def}
Функции образования $\xi$-множеств обозначаются как 

$\mathbf {\cap \left( {A,B} \right)}$ и $\mathbf {\cup \left( {A,B} \right)}$ 

\noindent
и имеют следующий вид:
%%%%%%%%%%%%%%%%%%%%%%%%%%%%%%%%%%%
\begin{equation}\label{cap}
\mathbf {\cap \left( {A,B} \right) = A \cap B \cup A \cap B \cup A \cap B \cup A \cap B \cup A \cap B \cup ... }
\end{equation}
\begin{equation}\label{cup}
\mathbf {\cup \left( {A,B} \right) = A \cup B \cap A \cup B \cap A \cup B \cap A \cup B \cap A \cup B \cap ...}
\end{equation}
%%%%%%%%%%%%%%%%%%%%%%%%%%%%%%%%%%%
\end{Def}

Используя теорему \ref{theorem2}, значения этих функций будут следующими $\xi$- множествами: 

$\mathbf { \cup \left( {A,B} \right) = \left( {A \cup B} \right)||A}$ 

\noindent
и 

$\mathbf {\cap \left( {A,B} \right) = A||\left( {A \cap B} \right)}$ соответственно. Из определения функций (\ref{cap}) и (\ref{cup}) следует следующая теорема для $\xi$-множества класса 2

\begin{Th}
Для двух множеств $\mathbf A$ и $\mathbf B$ которым равняется $\xi$ множество класса 2 справедливо утверждения: $\mathbf {A \subseteq B}$ или $\mathbf {B \subseteq A}$.
\end{Th}

\noindent
\textbf{Операции над $\xi$-множествами класса 2, образования $\xi$-множеств класса больше 2} 

Объединения $\xi$-множеств класса 2. Возьмем два $\xi$ множества класса 2 $\mathbf {A||F}$, $\mathbf {B||C}$ и объединим их $\mathbf {A||F \cup B||C}$. Так как каждое из них равняется одновременно $\mathbf {A,~F}$ и $\mathbf {B,~C}$ соответственно, то определим их объединение при условии 
%%%%%%%%%%%%%%%%%%%%%%%%%%%%%%%%%%%
$$\mathbf {A||F \cup B||C = \left\{ {\begin{array}{*{20}{c}}
	{\mathbf {A||F \cup B = \left\{ {\begin{array}{*{20}{c}}
			{\mathbf {A \cup B}} \\ 
			{\mathbf {F \cup B}} 
			\end{array}} \right.}} \\ 
	{\mathbf {A||F \cup C = \left\{ {\begin{array}{*{20}{c}}
			{\mathbf {A \cup C}} \\ 
			{\mathbf {F \cup C}} 
			\end{array}} \right.}} 
	\end{array}} \right.}$$

%%%%%%%%%%%%%%%%%%%%%%%%%%%%%%%%%%%
\noindent
что будет равняться $\xi$-множеству класса 4 
%%%%%%%%%%%%%%%%%%%%%%%%%%%%%%%%%%%%%%%%
$$\mathbf {A||F \cup B||C = \left( {A \cup B} \right)||\left( {B \cup F} \right)||\left( {A \cup C} \right)||\left( {F \cup C} \right)}$$

\noindent
пересечения и разность определяется по той же схеме
%%%%%%%%%%%%%%%%%%%%%%%%%%%%%%%%%%%
$$\mathbf {A||F \cap B||C = \left( {A \cap B} \right)||\left( {B \cap F} \right)||\left( {A \cap C} \right)||\left( {F \cap C} \right)}$$
%%%%%%%%%%%%%%%%%%%%%%%%%%%%%%%%%%%
$$\mathbf {A||F\backslash B||C = \left( {A\backslash B} \right)||\left( {B\backslash F} \right)||\left( {A\backslash C} \right)||\left( {F\backslash C} \right)}$$

Используя эти определения операций над $\xi$-множествами, докажем следующую теорему
\begin{Th}\label{theorem3}
При операциях объединения, пересечения и разность двух $\xi$-множеств классов $n$ и $m$, результирующее $\xi$-множество будет иметь класс $k \leqslant n m$.	
\end{Th}

\begin{proof}
	В самом деле, если рассмотреть определения операций объединения, пересечения и разность двух $\xi$-множеств классов 2 и взять к сведению что $\mathbf {A = A||A||...}$, то теорема доказана для этих $\xi$-множеств. Для $\xi$-множеств класса больше 2 теорема доказывается по индукции.
\end{proof}

Тем самым мы получили средство образования $\xi$-множества произвольного класса используя эти операции.

\noindent
\textbf{Принадлежность элементов к $\xi$-множеству} 

Возьмем произвольное $\xi$-множество 

\begin{equation}\label{A`}
{ \mathbf{{\tilde A}^\aleph }} = \mathop {||}\limits_{i \in \mathbf R} {\mathbf{A}_i}
\end{equation}

\noindent
класса $\mathbf{{\aleph} = \overline{\overline R} }$. Рассмотрим $\xi$-множество (\ref{A`}) класса $\aleph_0$ т.е. $\mathbf{\overline{\overline R}  = {\aleph _0}}$, принадлежность элементов к этому множеству неоднозначно т.е. если элемент $a$ принадлежит множествам $\mathbf{A}_i$ некоторого подмножества множества всех $\mathbf{A}_i$ из (\ref{A`}), то эта принадлежность обозначается как 
%%%%%%%%%%%%%%%%%%%%%%%%%%%%%%%%%%%
$$\mathbf{a\overset{{{k}_{1}},{{k}_{2}},...{{k}_{d}...}}{\mathop{\in }}\, \widetilde A^{{\aleph _0}} = {{A}_{1}}||{{A}_{2}}||...||{{A}_{n}}||...\Leftrightarrow a\in {{A}_{{{k}_{1}}}},a\in {{A}_{{{k}_{2}}}},...a\in {{A}_{{{k}_{d}}}} }$$ 

\noindent
или      
%%%%%%%%%%%%%%%%%%%%%%%%%%%%%%%%%%%
$$\mathbf{ a\mathop  \in \limits^{{k_1},{k_2},...{k_d}} {\kern 1pt} \widetilde A^{{\aleph _0}}}$$

\noindent
где индексы $\mathbf{{k}_{1},{k}_{2},...{k}_{d}...}$ над знаком принадлежности означают, каким именно множествам $\mathbf{A}_i$ из (\ref{A`}) принадлежит элемент $\mathbf a$.

Если мощность множества множеств $\mathbf{A}_i$ произвольна $\aleph$ т.е. $\xi$-множество имеет класс $\aleph$, то принадлежность элемента $\mathbf a$ к $\mathbf{{\widetilde A^\aleph }}$ определяется как 
%%%%%%%%%%%%%%%%%%%%%%%%%%%%%%%%%%%
$$\mathbf{a\mathop  \in \limits^T {\widetilde A^\aleph }}$$

\noindent
где $\mathbf{T = \left\{ {i|a \in {A_i}} \right\},~T \subset R,~\overline{\overline R} = \aleph}$.

Если элемент $b$ принадлежит всем множествам $\mathbf{A_i}$ из (\ref{A`}), то это обозначается как $\mathbf{b\mathop  \in \limits^{all} {\widetilde A^n}}$, если $c$ не принадлежит всем $\mathbf{A_i}$ из (\ref{A`}), то это обозначается как $\mathbf{c\mathop  \in \limits^0 {\widetilde A^n}}$ или в традиционном смысле как $\mathbf{c \notin {\widetilde A^n}}$.

\noindent
\textbf{Интерпретация $\xi$-множеств}

Так как рациональные числа есть что-то что находится межу целыми числами то $\xi$-множества это то что находится между обычными множествами. Обычное множество это то что находится <<между собой>>.

%%%%%%%%%%%%%%%%%%%%%%%%%%%%%%%%%%%%%%%%
%%%%%%%%%%%%%%%%%%%%%%%%%%%%%%%%%%%%%%%%
\noindent
\textbf{Аксиома $\xi$-множества, интерпретация парадокса Рассела} 

Добавим аксиому существования $\xi$-множества:
%%%%%%%%%%%%%%%%%%%%%%%%%%%%%%%%%%%%%%%%
$$\mathbf{\forall a \forall b  \exists \tilde a  (b \subseteq a \leftrightarrow   \tilde a = a \wedge \tilde a=b)}$$

\noindent
к системе аксиом наивной теории множеств и используем ее для интерпретации парадокса Рассела. Рассмотрим множества:

$\mathbf {U = \{X \mid X = X\}}$~--- множество всех множеств; 

$\mathbf{U_R = \{ X \mid X \notin X\}}$~--- множество Рассела, $\mathbf{{U_R} \in {U_R} \Leftrightarrow {U_R} \notin {U_R}}$;

$\mathbf{U_D = \{ X \mid X \in X\}}$;

$\mathbf{W = \{ X \mid X\mathop  \in \limits^1 X\}}$;

Примером множества принадлежащего $\mathbf W$ есть множество:   
%%%%%%%%%%%%%%%%%%%%%%%%%%%%%%%%%%%
$$\mathbf{{A_a} = \{ a\}  \cup \{ X \mid a \in X \wedge X \notin X\} }$$

\noindent
в самом деле, так как $\mathbf{a \in A_a}$, то справедливо утверждения 
%%%%%%%%%%%%%%%%%%%%%%%%%%%%%%%%%%%
$$\mathbf{{\bf{A_a}} \in {\bf{A_a}} \Leftrightarrow {\bf{A_a}} \notin {\bf{A_a}}}$$

Учитывая все свойства множеств  $\mathbf{U_R}$, $\mathbf W$ и $\mathbf{U_D}$ можно сделать следующие утверждения:

\begin{Th}\label{W}
	$\mathbf{\forall K \in {U_R}\left( {K \ne {U_R} \Leftrightarrow K \notin W} \right)}$
\end{Th}

\begin{proof}
	В самом деле для всех множеств принадлежащих к $\mathbf{U_R}$, кроме самого $\mathbf{U_R}$ (так как $\mathbf{{U_R}\mathop  \in \limits^1 {U_R}}$), справедливо утверждения $\mathbf{X \notin X}$, что не соответствует свойству определения множества $\mathbf W$.
\end{proof}

Используя теорему \ref{W}, представим множество Рассела:
%%%%%%%%%%%%%%%%%%%%%%%%%%%%%%%%%%%
$$\mathbf{{U_R} = A||B = \left( {U_R \cup \left( {{U_R}\backslash \left[ {W \cap {U_R}} \right]} \right)} \right)||\left( {{U_R}\backslash \left[ {W \cap {U_R}} \right]} \right)}$$

\noindent
\textbf{Гипотеза о равномощности множества и его булеана, Парадокс Кантора}

Парадокс Кантора \cite{11}~--- парадокс теории множеств, который демонстрирует, что предположение о существовании множества всех множеств ведет к противоречиям и, следовательно, противоречивой является теория, в которой построение такого множества возможно.

Для объяснение парадокса Кантора рассмотрим теорему кантора.

\begin{Th}[Кантор]\label{kantor}
Любое множество менее мощно, чем множество всех его подмножеств.
\end{Th}
\begin{proof}
Предположим, что существует множество $\mathbf A$, равномощное множеству всех своих подмножеств $2^\mathbf{A}$, то есть, что существует такая биекция $\mathbf f$, ставящая в соответствие каждому элементу множества $\mathbf A$ некоторое подмножество множества $\mathbf A$.

Рассмотрим множество $\mathbf B$, состоящее из всех элементов $\mathbf A$, не принадлежащих своим образам при отображении $\mathbf f$ (оно существует по аксиоме выделения): 

$$\mathbf{B=\left\{\,x\in A \mid x\not\in f(x)\,\right\}}$$ 

\noindent
$\mathbf f$ биективно, а $\mathbf{B \subseteq A}$, поэтому существует $\mathbf{y \in A}$ такой, что $\mathbf{f(y) = B}$. Теперь посмотрим, может ли $\mathbf{y}$ принадлежать $\mathbf B$. Если $\mathbf{y \in B}$, то $\mathbf{y \in f(y)}$, а тогда, по определению $\mathbf B$, $\mathbf{y \not\in B}$. И наоборот, если $\mathbf{y \not\in B}$, то $\mathbf{y \not\in f(y)}$, а следовательно, $\mathbf{y \in B}$. В любом случае, получаем противоречие.

Следовательно, исходное предположение ложно и $\mathbf A$ не равномощно $\mathbf{2^A}$.
\end{proof}

Доказательство теоремы \ref{kantor} производится методом от противного и на основе противоречия (с точки зрения наивной теории множеств) ${\bf{y}} \in {\bf{B}} \Leftrightarrow {\bf{y}} \notin {\bf{B}}$. Но если принять чо множество $\bf B$ это $\xi$-множество класса 2:
%%%%%%%%%%%%%%%%%%%%%%%%%%%%%%%%%%%%%%%%
$$\mathbf{B = \{ y\}  \cup D||D} \text{, где}\mathbf{D = \{ x \mid x \notin f(x) \wedge x \ne y\} }$$

\noindent
тогда $\mathbf{y\mathop  \in \limits^1 B}$.

Гипотеза о равномощности множества и его булеана заключается в следующем:

Множество всех множеств не может существовать в терминах наивной теории множеств, но если рассматривать его с точки зрения $\xi$-множеств, то имеет место равномощности множества и его булеана.

\end{document}